\documentclass[a4paper]{article}

\usepackage{theorem}
\usepackage{latexsym}
\usepackage{amsmath}
\usepackage{amssymb}
\usepackage{graphicx}
\graphicspath{{figures/}}
\usepackage{gastex}
\usepackage[ruled,vlined]{algorithm2e}

{
\theorembodyfont{\rmfamily}
\newtheorem{example}{Example}[section]
\newtheorem{note}[example]{Note}
}

\newtheorem{theoreme}[example]{Theorem}
\newtheorem{corollaire}[example]{Corollary}

\newtheorem{definition}[example]{Definition}
\newtheorem{proposition}[example]{Proposition}
\newtheorem{lemme}[example]{Lemma}
\newcommand{\qed}{\hfill$\Box$ \vspace{0.5 cm}}
\newenvironment{proof}{{\it Proof~: }}{\qed}

\newcommand{\liste}[2][n]{\ensuremath{#2_1,\ldots,#2_{#1}}}
\newcommand{\shot}{shot-set}
\newcommand{\sh}[1]{\ensuremath{s(#1)}}
\newcommand{\mylabel}[1]{\label{#1}}
\newcommand{\ie}{\emph{i.e.}}

\newcommand{\lcfg}{\ensuremath{L(CFG)}}
\newcommand{\lasm}{\ensuremath{L(ASM)}}
\newcommand{\lmcfg}{\ensuremath{L(MCFG)}}
\newcommand{\confspace}[1]{\ensuremath{L(#1)}}
\newcommand{\ens}[1]{\ensuremath{\{#1\}}}

\newcommand{\property}[1]{Proposition~\ref{#1}}

\newcommand{\gr}{Grounding Modification}
\newcommand{\mul}{Multiplying Modification}

\newcommand{\dbot}[1]{\ensuremath{d_{\bot}(#1)}}
\newcommand{\dmoins}[2][]{\ensuremath{d_{#1}{^-}(#2)}}
\newcommand{\dplus}[2][]{\ensuremath{d_{#1}^+(#2)}}
\newcommand{\dinf}[2][]{\ensuremath{d_{#1}{^<}(#2)}}
\newcommand{\dsup}[2][]{\ensuremath{{d_{#1}}^>(#2)}}
\newcommand{\setedge}[2]{\ensuremath{S_{#1}^{#2}}}

\newcommand{\confinit}[1]{\ensuremath{\sigma_0(#1)}}
\newcommand{\conffin}[1]{\ensuremath{\sigma_f(#1)}}

\def\N{\mbox{I\hspace{-.15em}N}}

\title{Classes of Lattices Induced by Chip Firing \\(and Sandpile) Dynamics}
\author{Cl\'emence Magnien}
\date{}

\begin{document}
\maketitle

\begin{abstract}
In this paper we study
three classes of models widely used in physics, computer science and social science: the
Chip Firing Game,
the Abelian Sandpile Model and the Chip Firing Game on a mutating graph.
We study the set of configurations reachable from a given initial configuration, called the \emph{configuration space} of a model,
and try to determine the main properties of such sets.
We study the order induced over the configurations by the evolution rule.
This makes it possible to
 compare the power of expression of these models.
It is known that the configuration spaces we obtain are lattices, a special kind of partially ordered set.
Although the Chip Firing Game on a mutating graph is a generalization of the usual Chip Firing Game, we prove that these models 
generate exactly the same configuration spaces.
We also prove that
the class
of lattices induced by the Abelian Sandpile Model is strictly included
in the class of lattices induced by the Chip Firing Game, but contains
the class of distributive lattices, a very well known class.
\end{abstract}

\section*{Introduction}
The Chip Firing Game (CFG), the Abelian Sandpile Model (ASM), and the
Chip Firing Game on a mutating graph (which we will call Mutating Chip
Firing Game or MCFG) are closely related models studied in physics~\cite{BTW87,DRSV95}, computer science~\cite{GMP98a,BLS91,BL92} and social science~\cite{Big97,Big99,Heu99}.
These three models are variations of the following game:
given a graph and a distribution of chips on its vertices
(called \emph{configuration}), one may select a vertex that contains at least as many chips as its outdegree, and move one chip from this vertex 
along each of its outgoing edges, to the vertex at the other extremity.
Many questions have naturally arisen and given matter for research about this game:
given a graph and an initial configuration, does the game stop after some time, or can it be played forever~\cite{BLS91,BL92,Eri96}~?
For a given graph, what are the properties of the configurations such that no move is possible~? This has led to the algebraic study of some of these configurations called \emph{recurrent configurations}~\cite{DRSV95,CR00,Big99}.
Given a graph, what can be said about the set of the configurations that can be reached from a given initial configuration~\cite{LP00,MPV01}~?

This paper is devoted to the study of this last question.
 Our purpose is to study the
set of  configurations reachable from a given initial configuration, 
 which we will call
the \emph{configuration space} of the game.
This set is naturally ordered by the relation of reachability induced by the evolution rule.
In \cite{MPV01} it was proved that the configuration spaces induced by CFGs have a very strong structure, and that the class they form is situated between two very well-known classes of orders
(the distributive and the ULD lattices).
Our aim is to determine if the same kind of results also holds for the two other models.
We will try to determine the differences between the classes of configuration spaces induced by each model, \ie{}
given the configuration
space of a game, we will try to decide if it is isomorphic to the
configuration space of a game of another type.
For instance, given the configuration space of a CFG, does there exist
 a MCFG or
an ASM such that its configuration space is isomorphic to it~?

In Section \ref{secrecalls}, we give the definitions and known results used in this paper, as well as the definitions of the models.
In Section \ref{seccfgmut}, we show that CFGs and MCFGs induce exactly  the same class of configuration spaces,
and in Section \ref{seccfgasm}, we will 
compare the configuration spaces of CFGs and ASMs.
Since an ASM is a special CFG, we will 
study under which condition a CFG can be transformed into an ASM.
We will give a sufficient (but not necessary) condition of such a transformation, and
we will give an example of a CFG that cannot be transformed into an ASM.

\section{Definitions and known results} \mylabel{secrecalls}

We will need some definitions from order and lattice theory to describe the configuration spaces of the models.
We give them first, after what we give the definition of each model,
 as well as the previously known results about them.
\subsection{Posets and lattices}

A \emph{partially ordered set} (or \emph{poset}) is a set equipped with
an order relation $\le$ (\ie{} a transitive, reflexive and
antisymmetric relation).  If $x$ and $y$ are two elements of a poset, we say
that $x$ is \emph{covered} by $y$ (or $y$ \emph{covers} $x$), and
write $x\prec y$ (or $y\succ x$) if $x<y$ and $x\le z<y$ implies $z=x$.
We then say that $x$ is a \emph{lower cover} of $y$ (or that $y$ is an \emph{upper cover} of $x$).
The \emph{interval} $[x,y]$ is the set \ens{z, x\le z\le y}.  To
represent a poset $P$ we will use its Hasse diagram, defined as
follows~:
\begin{itemize}
\item each element $x$ of $P$ is represented by a point $p_x$ of the plane,
\item if $x<y$, then $p_x$ is lower than $p_y$, and
\item $p_x$ and $p_y$ are joined by a line if and only if $x\prec y$.
\end{itemize}
Two posets $P$ and $P'$ are \emph{isomorphic} if there exists a bijection $\varphi: P \longrightarrow P'$
satisfying: for all $x,y\in P$, $x\le y \iff \varphi(x) \le \varphi(y)$.

A poset $L$ is a \emph{lattice} if any two elements $x,y$ of $L$ have
a least upper bound (called \emph{join} and denoted by
$x\vee y$) and a greatest lower bound
(called \emph{meet} and denoted by $x\wedge y$).
The join
$x\vee y$ is the (unique) smallest element greater than both $x$ and $y$, and
$x\wedge y$ is defined dually.
All
the lattices considered here are finite, therefore they have a unique maximal and a unique minimal element.

A lattice is a \emph{hypercube of dimension $n$} if it is isomorphic
to the set of all subsets of a set of $n$ elements, ordered by
inclusion.  
Hypercubes are also called \emph{boolean lattices}.
A lattice is \emph{upper locally distributive} (denoted by
\emph{ULD}) \cite{Mon90} if the interval between any element and the
join of all its upper covers is a hypercube.
\emph{Lower locally distributive} (LLD) lattices are defined dually.
All ULD lattices are
ranked,
\ie{} all the paths in the covering relation from the minimal to the maximal element have the same length.
A lattice $L$ is \emph{distributive} if it 
satisfies one of the two following relations of distributivity (that are equivalent and imply each other):
$$\mbox{for all }x,y,z\in L, x\vee(y\wedge z)=(x\vee y)\wedge(x\vee z)$$
$$\mbox{for all }x,y,z\in L, x\wedge(y\vee z)=(x\wedge y)\vee(x\wedge z).$$
A distributive lattice is a lattice that
is at the same time
upper \emph{and} lower locally distributive, \ie{} if 
the interval between any element and, on the one hand the join of its upper covers, and on the other hand
 the meet of all its lower
covers both are hypercubes.
In the sequel we will only be concerned with distributive and ULD lattices.
Figure~\ref{figexordres} shows examples of posets and different types of lattices.

\begin{figure}
\begin{center}
\includegraphics{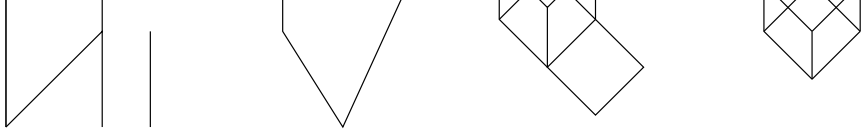}
\end{center}
\caption{From left to right: a poset, a lattice, a distributive lattice and a hypercube of dimension 3}
\mylabel{figexordres}
\end{figure}

For a more complete introduction to posets and lattices, see for in\-stance~\cite{DP90}.

\subsection{The different kinds of Chip Firing Games}
In this section we give the definitions of the Abelian Sandpile Model, the Chip Firing Game and the Mutating Chip Firing Game.
We begin by presenting the features shared by the three models.
Then we detail what is specific about each of them.
We will see that some models are generalizations of others.
Finally, we give some results about the CFG which will be useful in this paper.

\subsubsection{Definitions}
\mylabel{subsubdef}
Each model is defined over a  graph $G=(V,E)$, called the \emph{support graph} of the game
(undirected graphs will be regarded as directed by replacing each undirected edge \ens{i,j} by the two directed edges $(i,j)$ and $(j,i)$).
All graphs  are supposed to be multigraphs, \ie{} 
multiple edges between two vertices are allowed, 
therefore all edge sets are also supposed to be multisets.
A \emph{configuration} of the game is 
a mapping $\sigma:V\mapsto \N$ which associates
 a weight to each vertex;
this weight may be considered as a number of \emph{chips} stored in the vertex.
The game is played with respect to
 the following evolution rule, also called the \emph{firing} rule:
if a vertex $v$ contains at least as many chips as its outdegree,
we can transfer a chip from $v$ along each of its outgoing edges to the corresponding vertex.
We call this process \emph{firing $v$}.
 If $\sigma$ is the configuration we start from, and $\sigma'$ is obtained from $\sigma$ by firing $v$,
we write $\sigma\stackrel{v}{\longrightarrow} \sigma'$, and we call $\sigma$ a \emph{predecessor} of $\sigma'$.

\begin{note}
We consider that the firing rule cannot be applied to 
a vertex with no outgoing edges
 because 
firing it does not change the
configuration of the game and is therefore of no interest to us.
\end{note}

The Chip Firing Game (CFG) \cite{BL92} is defined over a directed graph.
We give an example of a CFG together with its configuration space in 
Figure~\ref{figcfgtreillis}.
\begin{figure}
\mbox{} \hfill
\begin{minipage}{3cm}
\scalebox{0.9}{
\begin{picture}(30,40)
\node(a)(0,0){1}
\node(b)(30,0){1}
\node(c)(20,20){1}
\node(d)(5,30){0}
\gasset{ExtNL=y,NLdist=1}
\nodelabel[NLangle=180](a){$a$}
\nodelabel[NLangle=0](b){$b$}
\nodelabel[NLangle=0](c){$c$}
\nodelabel[NLangle=180](d){$d$}
\drawedge(a,c){}
\drawedge(b,c){}
\drawedge[curvedepth=2](c,d){}
\drawedge[curvedepth=-2](c,d){}
\end{picture}}
\end{minipage}
\hfill
\begin{minipage}{3.5cm}
\scalebox{0.6}{\input{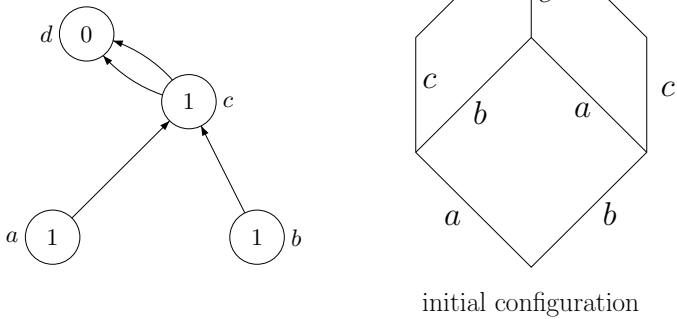}} 
\end{minipage}
\hfill \mbox{}
\caption{An initial configuration of a CFG  and its configuration space. Each edge $(\sigma,\sigma')$ is labelled by the name of the vertex fired to reach $\sigma'$ from $\sigma$.}
\mylabel{figcfgtreillis}
\end{figure}
The Abelian Sandpile Model (ASM) \cite{BTW87} is defined over an undirected graph,
with a distinguished vertex called the \emph{sink} and denoted by $\bot$. 
The sink can never be fired.
Finally the Chip Firing Game on a mutating graph (MCFG) \cite{Eri96} is played on a directed graph that changes during the game in the following way: after the firing of a vertex, its outgoing edges are removed and new ones are added
in a pre-determined way
(the vertex set $V$ remains the same). 
A \emph{position} of a MCFG consists in: 
\begin{itemize}
\item a directed graph on $V$,
\item on each node $v\in V$ a nonnegative number of chips,
\item for each node $v\in V$ an infinite sequence $M_v$ of multisets of nodes in $V$, which is the mutation sequence of $v$.
\end{itemize}
A vertex $v$ may be fired if the number of chips on $v$ is at least \dplus{v} (the current number of outgoing edges of $v$). 
The mutation of the graph takes place after the firing of the vertex.
The outgoing edges of $v$ are removed, and
if $M_v = (\setedge{v}{1},\setedge{v}{2},\ldots)$ is the mutation sequence of $v$,
then a new edge $(v,w)$ is created for each $w\in \setedge{v}{1}$.
Finally, \setedge{v}{1} is removed from $M_v$, so that $M_v$ becomes 
$(\setedge{v}{2},\setedge{v}{3},\ldots)$.
We define the \emph{initial support graph} of a MCFG to be its support graph in the initial position.
An example of a firing sequence in a MCFG is given in Figure~\ref{figcfgmut}.

\begin{figure}
\hfill
\begin{minipage}{2cm}
\scalebox{0.65}{\begin{picture}(35,35)
\node[fillgray=0.85](a)(5,5){1}
\node(b)(30,5){0}
\node(c)(17,30){3}
\gasset{ExtNL=y,NLdist=1}
\nodelabel[NLangle=180](a){$a$}
\nodelabel[NLangle=0](b){$b$}
\nodelabel[NLangle=90](c){$c$}
\drawedge[curvedepth=-2](a,b){}
\drawedge[curvedepth=2](b,c){}
\drawedge[curvedepth=2](c,b){}
\drawedge[curvedepth=-2](c,a){}
\end{picture}}
\end{minipage}
\hfill 
\begin{minipage}{0.5cm}
$\longrightarrow$ 
\end{minipage}
\hfill
\begin{minipage}{2cm}
\scalebox{0.65}{\begin{picture}(35,35)
\node(a)(5,5){0}
\node[fillgray=0.85](b)(30,5){1}
\node(c)(17,30){3}
\gasset{ExtNL=y,NLdist=1}
\nodelabel[NLangle=180](a){$a$}
\nodelabel[NLangle=0](b){$b$}
\nodelabel[NLangle=90](c){$c$}
\drawedge[curvedepth=-2](a,c){}
\drawedge[curvedepth=2](b,c){}
\drawedge[curvedepth=2](c,b){}
\drawedge[curvedepth=-2](c,a){}
\end{picture}}
\end{minipage}
\hfill 
\begin{minipage}{0.5cm}
$\longrightarrow$ 
\end{minipage}
\hfill
\begin{minipage}{2cm}
\scalebox{0.65}{\begin{picture}(35,35)
\node(a)(5,5){0}
\node(b)(30,5){0}
\node[fillgray=0.85](c)(17,30){4}
\gasset{ExtNL=y,NLdist=1}
\nodelabel[NLangle=180](a){$a$}
\nodelabel[NLangle=0](b){$b$}
\nodelabel[NLangle=90](c){$c$}
\drawedge[curvedepth=-2](a,c){}
\drawedge[curvedepth=2](b,a){}
\drawedge[curvedepth=2](c,b){}
\drawedge[curvedepth=-2](c,a){}
\end{picture}}
\end{minipage}
\hfill 
\begin{minipage}{0.5cm}
$\longrightarrow$ 
\end{minipage}
\hfill
\begin{minipage}{2cm}
\scalebox{0.65}{\begin{picture}(35,35)
\node(a)(5,5){1}
\node(b)(30,5){1}
\node(c)(17,30){2}
\gasset{ExtNL=y,NLdist=1}
\nodelabel[NLangle=180](a){$a$}
\nodelabel[NLangle=0](b){$b$}
\nodelabel[NLangle=90](c){$c$}
\drawedge[curvedepth=-2](a,c){}
\drawedge[curvedepth=2](b,a){}
\drawedge[curvedepth=2](c,b){}
\drawedge[curvedepth=-2](c,a){}
\end{picture}}
\end{minipage}
\hfill \mbox{}
\caption{An example of a firing sequence in a MCFG, where the beginning of the mutation sequences of the vertices are : 
$M_a=(\ens{c},\ldots),M_b=(\ens{a},\ldots),M_c=(\ens{a,b},\ldots)$}
\mylabel{figcfgmut}
\end{figure}
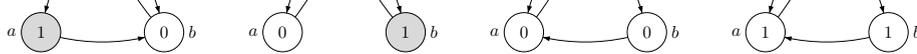

These three games are strongly convergent games \cite{Eri93}, which implies that, given an initial configuration, 
either a given game can be played forever
or it reaches a unique fixed point (where no firing is possible),
called the \emph{final configuration},
which does not depend on the order in which the vertices are fired. 
We will only consider \emph{convergent} games,\ie{} games that reach a fixed point.
We call \emph{execution} of a game any sequence of firing that, 
starting from the initial configuration, reaches the final configuration.

We call \emph{configuration space} of a game $C$,
and we denote it by \confspace{C}, the set of all the configurations reachable from the initial configuration, ordered by the reflexive and transitive closure of
the predecessor relation.
It is known \cite{BLS91,LP00,Eri96} that when a game is convergent,
its configuration space
is a ULD lattice.
This is a very strong property, because ULD lattices are very structured sets.
For instance, an immediate consequence of this is the fact that all the firing sequences, in any convergent ASM, CFG or MCFG, from the initial configuration to the final configuration, have the same length.
We will say that two convergent
games $C$ and $C'$ are
\emph{equivalent} if \confspace{C} is isomorphic to \confspace{C'}.

We denote by \lcfg{}, \lasm{} and \lmcfg{}
the classes of lattices that are the configurations spaces of convergent CFGs, ASMs and MCFGs respectively.
  It is possible to
guarantee that a CFG is convergent by the presence in the support
graph of a sink (a vertex with no outgoing edges) reachable from all vertices
\cite{LP00}.
This also holds for MCFGs.
In the ASM, the presence of the distinguished vertex  which can never be fired
(called the sink
because it has the same function as in the directed case)
 guarantees that the game is always convergent.

These three models are very close to one another, and it is easy to
see that some of these are generalizations of others (the class of
lattices they induce are included in one another):  since undirected
graphs are particular directed graphs, one can consider an ASM as a
particular CFG, where each undirected edge $\ens{u,v}$ with $u,v\not=\bot$ is replaced by two opposite
directed edges $(u,v)$ and $(v,u)$, 
and where each undirected edge \ens{v,\bot} is replaced only by the directed edge $(v,\bot)$.
 Therefore, we
obtain that $\lasm\subseteq\lcfg$.  
Likewise, a CFG can be regarded as a MCFG where the graph remains the same after each mutation, therefore $\lcfg\subseteq\lmcfg$.
Our aim is to study the other relations between these classes of lattices.

\subsubsection{Previous results}
Now let us give some definitions and known results about CFGs, useful for simplifying the notations and proofs in the sequel.
\begin{definition}
A convergent ASM, CFG or MCFG is \emph{simple} if, during an execution, each vertex is fired at
most once.
\end{definition}

\begin{theoreme} \cite{MPV01} \mylabel{thmsimple}
Any convergent CFG is equivalent to a simple CFG.
\end{theoreme}

Thanks to this result, all the CFGs considered in the sequel  will be simple, and
such that their support graph has exactly one sink (denoted by $\bot$)
and such that all vertices except $\bot$ are fired during an execution.
This is always possible because, if this is not the case:
\begin{itemize}
\item either there is no sink and all vertices are fired during an execution; then we can add an isolated vertex to the graph, which becomes the sink,
\item or there exists a vertex $v$ that is never fired during any execution but is not a sink; then we can remove its outgoing edges without changing the configuration space of the CFG, and $v$ becomes the sink, 
\item or there is more than one sink; then we can merge them into a single vertex without changing the configuration space.
\end{itemize}
In the sequel, we will restrict ourselves to CFGs for which
 the support graph has no loops (\ie{} no $(v,v)$ edges).
This is always possible because, for a simple CFG $C$, if there is a loop on a vertex $v$
($v$ cannot be the sink), we can replace it by an edge $(v,\bot)$,
and the resulting CFG is equivalent to $C$, since $C$ is \emph{simple}.

Given a vertex $v$, we denote by $d^-(v)$ its indegree, by \dplus{v} its outdegree, by \dbot{v} the number of edges from $v$ to $\bot$, and we define $d(v)=d^+(v)-\dbot{v}$.
Given two vertices $u$ and $v$,
we denote by $d(u,v)$ the number of edges from $u$ to $v$.
Given a CFG and a vertex $v$ of its support graph,
the initial number of chips in $v$ is denoted by \confinit{v}, and the number of chips in $v$ in the final configuration by \conffin{v}
(since all the CFGs considered are simple, we have $\conffin{v}=\confinit{v}+d^-(v)-d^+(v)$).
The number of chips \emph{needed to fire $v$} is the difference
between \confinit{v} and \dplus{v}, \ie{} 0 if $\confinit{v}\ge\dplus{v}$,
and $\dplus{v} - \confinit{v}$ otherwise.

The configuration spaces of simple CFGs can be described more easily than in the general case. Indeed we have
the following results:
\begin{lemme} \cite{MPV01,LP00}
In a simple CFG, if, starting from the same configuration, two sequences of firing 
lead to the same configuration, then the vertices fired in each sequence are the same.
\end{lemme}
This allows us to define the \emph{\shot{}} \sh{\sigma} of a configuration $\sigma$ as the set of vertices fired to reach $\sigma$ from the initial configuration.
Given a CFG with support graph $(V,E)$, we say that a subset $X\subseteq V$ is a \emph{valid} \shot{} if there exists a configuration $\sigma$ reachable from the initial configuration
such that $\sh{\sigma}=X$.
A list (\liste{v}) of vertices is a \emph{valid firing sequence} if, for each $i$, \ens{\liste[i]{v}} is a valid \shot{}.
The configuration space of a CFG is isomorphic to the lattice of the \shot{}s of its configurations, ordered by inclusion.
The join of any two elements $a$ and $b$ is given by the following formula \cite{LP00}:
$$\sh{a\vee b}=\sh{a}\cup\sh{b}.$$

Here we give another way to characterize the configuration space of a CFG with respect to the \shot{}s:
given a CFG and its vertex set, we can associate
with each vertex $v$,  the set of the configurations in which $v$ can be fired.
Among these, we distinguish the smallest configurations (the ones such that their \shot{} is minimal with respect to the inclusion), and we say that 
their \shot{} represent the \emph{first times} at which $v$ can be fired.
For instance, in Figure \ref{figcfgtreillis}, 
the \shot{}s of the minimal configurations in which $c$ can be fired are $\{a\}$ and $\{b\}$,
 and we say that the first times at which $c$ can be fired are $a$ and $b$.
The minimal configuration in which $a$ and $b$ can be fired is the minimal point of the lattice, and its \shot{} is the empty set.
We say that the first time at which $a$ and $b$ can be fired is the beginning of any execution.
The knowledge, for all vertices, of the first times at which they can be fired, is a characterization of the configuration space of a CFG, as it is stated in the following proposition:

\begin{proposition}\mylabel{propcaracconfspace}
Let $C$ be a simple CFG and let $L=\confspace{C}$.
Let \ens{\liste{v},\bot} be the set of vertices of $C$.
For each $i,\in\ens{1,\ldots,n}$, let $X_i$ be the set representing the first times at which $v_i$ can be fired
(\ie{} $X_i$ is the set of the \shot{}s of the minimal configurations at which $v_i$ can be fired).
Then $L$ is completely determined by \ens{\liste{v}} and \ens{\liste{X}}.
\end{proposition}
\begin{proof}
We will construct $L$ from \ens{\liste{v}} and \ens{\liste{X}}.
Notice that, for each $i=0,\ldots,n$, and for each $x\in X_i$,
the set $\ens{v_i}\cup x$ is a valid \shot{} of $C$:
by definition, $x$ is the \shot{} of a configuration in which
$v_i$ can be fired.
We claim that all subsets $Y$ of \ens{\liste{v}} satisfying the following condition:
$$\mbox{for all } v_i \in Y, \mbox{ there exists } Z\in X_i,
\mbox{ such that } Z\subseteq Y$$
 are valid \shot{}s of elements of $L$.
Let $Y$ be such a set.
$Y$ can be decomposed into a union of sets \liste[k]{Y} such that,
for each $j=1,\ldots,k$,
there exists $i$ such that
$Y_j$ is the union of \ens{v_i} and of an element $x$ of $X_i$
(notice that the sets \liste[k]{Y} are not a partition of $Y$, some sets may overlap).
Therefore, each of the sets $Y_j$ is a valid \shot{} of $C$.
Since the \shot{}s are stable by union, $Y$ is itself a valid \shot{} of $C$.
On the other hand, every \shot{} of $C$ satisfies the condition above.

Therefore the set of subsets of \ens{\liste{v}} that satisfy this condition is the set of \shot{}s of $C$, which, ordered by inclusion, is isomorphic to $L$.
\end{proof}

This result will be useful in the sequel when we have to guarantee that a modification made to a game does not change its configuration space.
We now have all the definitions and tools needed to compare the classes of lattices induced by the models we have introduced.

\section{Comparison of MCFGs and CFGs} \mylabel{seccfgmut}
In this section we show that the CFGs and the MCFGs induce exactly the same class of configuration spaces.
We already know from Section \ref{subsubdef} that any CFG is equivalent to a MCFG, therefore 
we need to show that any MCFG is equivalent to a CFG.\\

We begin by proving that any MCFG is equivalent to a simple one.
If a MCFG is not simple, then there exists a vertex that is fired more than once during an execution.
In \cite{Eri96} it is shown that the number of times this vertex is fired is the same in \emph{every} possible execution of the game:
this vertex is fired the same number of times in all the firing sequences that lead from the initial to the final configuration.
We need this result to prove the
next theorem, which is an equivalent of Theorem~\ref{thmsimple} for MCFGs.
\begin{theoreme} \mylabel{cormutesimple}
Any convergent MCFG is equivalent to a \emph{simple} MCFG.
\end{theoreme}
\begin{proof}
Let $C$ be a MCFG with vertex set $V$ and
initial configuration $\sigma$.
Suppose $C$ is not simple: there exists a vertex $a$ 
 which is fired twice or more during any execution of $C$.
The idea of the proof is to split $a$ into two vertices $a_0$ and $a_1$ 
which will be fired alternatively (the first firing of $a$ in $C$ corresponds to a firing of $a_0$ in $C'$, the second to a firing of $a_1$, 
and so on),
so that each of them is fired strictly less often than $a$ during an execution.
Therefore, by iterating this process, one can transform $C$ into a simple MCFG.

Let $G$ be the initial support graph of $C$.
If for each vertex $v$, the sequence of mutation is 
$M_v=(\setedge{v}{1},\setedge{v}{2},\setedge{v}{3},\ldots)$, we define \setedge{v}{0} to be the set of outgoing edges of $v$ in $G$.
Therefore from now on, a MCFG will be defined by a set of vertices, an initial configuration, and  for each vertex $v$ 
an inifinite sequence of edge sets $M_v=(\setedge{v}{0},\setedge{v}{1},\setedge{v}{2},\ldots)$.
At any time in the game,
if a vertex $v$ with mutating sequence $(\setedge{v}{0},\setedge{v}{1},\setedge{v}{2},\ldots)$ has been fired	
$i$ times, then the current outgoing edges of $v$ are those given by 
$\setedge{v}{i}$.
We then
 denote
by $l_i(v)$ the number of loops on $v$ (\ie{} the number of occurrences of $v$ in \setedge{v}{i}),
by \dsup[i]{v} the number of edges going out of $v$ which are not loops (\emph{i.e.} $\dsup[i]{v}=\dplus[i]{v}-l_i(v)$),
and we define \dinf[i]{v} dually: $\dinf[i]{v}=\dmoins[i]{v}-l_i(v)$).
Finally we define $N$ to be twice the number of chips present in the game.

Let us define the MCFG $C'$ in the following way:
\begin{itemize}
\item the vertex set of $C'$ is $V'=V\backslash \{a\} \cup \{a_0,a_1\}$, with $a_0\not\in V$ and
  $a_1\not\in V$.
\item for each $v\not=a\in V$, and for each $i\in \N$ we define
   $\setedge{v}{i}$, the $i$-th set in the mutation sequence of $v$, in the
   following way:
\begin{itemize}
\item for each occurrence of a vertex $v'\not=a$ in $\setedge{v}{i}$, there are two occurrences of $v'$ in $\setedge{v}{i}$,
\item for each occurrence of $a$ in $\setedge{v}{i}$, there is one occurrence of $a_0$ and one occurrence of $a_1$ in $\setedge{v}{i}$.
\end{itemize}
\item for each $i\in\N$, we define the sets $\setedge{a_0}{i}$ and 
\setedge{a_1}{i},
the $i$-th sets in the mutating sequences of $a_0$ and $a_1$, in the following way:
\begin{itemize}
\item for each occurrence of a vertex $w\not=a$ in $\setedge{a}{2i}$ (resp. 
$\setedge{a}{2i+1}$), there are two occurrences of $w$ in $\setedge{a_0}{i}$
(resp. $\setedge{a_1}{i}$).
\item for each occurrence of $a$ in $\setedge{a}{2i}$ (resp. \setedge{a}{2i+1}), there is one occurrence of $a_0$ in  \setedge{a_0}{i} 
(resp. one occurrence of $a_1$ in  \setedge{a_1}{i}).
\item there are $N-\dsup[2i]{a}$ occurrences of $a_1$ in 
$\setedge{a_0}{i}$, and 
   $N-\dsup[2i+1]{a}$ occurrences of $a_0$ in $\setedge{a_1}{i}$.
\end{itemize}
\item for any vertex $v\not=a$, $\sigma'(v)=2\sigma(v)$.
\item $\sigma'(a_0)=\sigma(a)+N$.
\item $\sigma'(a_1)=\sigma(a)$.
\end{itemize}
An example of the construction of the graph is given in Figure \ref{figmcfgsimple}.

\begin{figure}
\begin{center}
\begin{tabular}{ccc}
$v$ in $C$ after $i$ firings & & $v$ in $C'$ after $i$ firings \\
\begin{minipage}{35mm}
\begin{center}
\begin{picture}(35,25)
\put(0,0){\framebox(35,25){}}
\node(v)(10,20){$v$} \node(v')(25,20){$v'$} \node(a)(18,5){$a$}
\drawedge(v,v'){} \drawedge(v,a){}
\end{picture}
\end{center}
\end{minipage}
& 
\begin{minipage}{1.5cm}
\begin{center}
$\longrightarrow$ \end{center}
\end{minipage}
&
\begin{minipage}{35mm}
\begin{center}
\begin{picture}(35,25)
\put(0,0){\framebox(35,25){}}
\node(v)(10,20){$v$} \node(v')(25,20){$v'$} 
\node(a0)(10,5){$a_0$} \node(a1)(25,5){$a_1$}
\drawedge[curvedepth=2](v,v'){} \drawedge[curvedepth=-2](v,v'){}
\drawedge(v,a0){} \drawedge(v,a1){}
\end{picture} 
\end{center}
\end{minipage} \\
\\
$a$ in $C$ after $2j$ firings & & $a_0$ in $C'$ after $j$ firings \\
\begin{minipage}{35mm}
\begin{center}
\begin{picture}(35,25)
\put(0,0){\framebox(35,25){}}
\node(a)(20,20){$a$} \node(w)(20,5){$w$}
\drawloop[loopangle=180](a){}
\drawedge(a,w){}
\end{picture}
\end{center}
\end{minipage}
& 
\begin{minipage}{1.5cm}
\begin{center} $\longrightarrow$ \end{center}
\end{minipage}
&
\begin{minipage}{35mm}
\begin{center}
\begin{picture}(35,25)
\put(0,0){\framebox(35,25){}}
\node(a0)(14,20){$a_0$} \node(a1)(29,20){$a_1$} \node(w)(20,5){$w$}
\drawloop[loopangle=180](a0){}
\drawedge[curvedepth=2](a0,w){} \drawedge[curvedepth=-2](a0,w){}
\drawedge(a0,a1){{\tiny $N-1$}}
\end{picture}
\end{center}
\end{minipage} \\
\end{tabular}
\end{center}
\caption{The outgoing edges of the vertices in the MCFG constructed during the simplification of a MCFG.}
\mylabel{figmcfgsimple}
\end{figure}
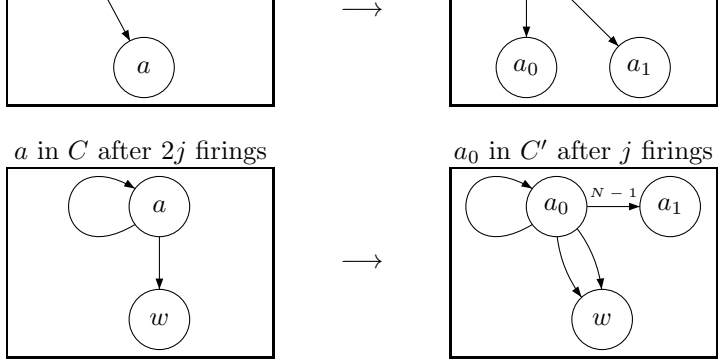

We call the \emph{double of $C$} the CFG obtained from $C$ by multiplying by two the initial configuration of all vertices, as well as the number of edges between each pair of vertices.
This CFG obviously is equivalent to $C$, and in the sequel we will show that $C'$ is equivalent to the double of $C$.
Notice that, except the two vertices $a_0$ and $a_1$, $C'$ is the
same as the double of $C$.
So if the firings of $a_0$ and $a_1$ take place in turn as
stated before, then we can easily show that
\confspace{C'} is isomorphic to \confspace{C}.

We first prove by induction that any firing sequence $s$ of $C$, leading to a configuration $\sigma$, corresponds to a valid firing sequence $s'$ in $C'$
(by replacing the occurrences of $a$ by occurrences of $a_0$ and $a_1$ in turn), leading to a configuration $\sigma'$, satisfying 
Proposition~\ref{prop1} and Proposition~\ref{prop2} or \ref{prop2}$'$.
\begin{gather}
\mbox{for all}\ v\neq a_0, a_1,\sigma'(v)=2\sigma(v) \mylabel{prop1}\\
\sigma'(a_0)=N+\sigma(a) \mbox{, and } \sigma'(a_1)=\sigma(a)  \mylabel{prop2}\\
\sigma'(a_0)=\sigma(a) \mbox{, and } \sigma'(a_1)=N+\sigma(a) \tag{\ref{prop2}$'$}
\end{gather}

The initial configuration of $C$ and $C'$ satisfy Propositions~\ref{prop1} and~\ref{prop2}.
So does any configuration reached by a firing sequence that does not contain $a$, 
because, except $a_0$ and $a_1$, $C'$ is the same as the double of $C$, and because
for each
$v\not= a_0,a_1$, there is the same number of edges from $v$ to $a_0$
than from $v$ to $a_1$.

Let now $s$ be a firing sequence of $C$ that does not contain $a$, but leading to a configuration $\sigma$ where $a$ can be fired.
$s$ corresponds to a firing sequence $s'$ in $C'$, leading to a configuration $\sigma'$ satisfying \property{prop1} and \property{prop2}.
In $\sigma'$ $a_0$ contains $N+\sigma(a)$ chips.
By construction, the outdegree of $a_0$ is:
$2\dsup[0]{a}+l_{0}(a)+(N-\dsup[0]{a}) = N + \dplus[0]{a}$.
Therefore, the fact that $a$ can be fired in configuration $\sigma$ implies that $a_0$ can be fired in configuration $\sigma'$.
Let $\sigma_2$ and ${\sigma_2}'$ be the configurations reached from $\sigma$ and $\sigma'$ by firing $a$ and $a_0$.
We will show that these configurations satisfy \property{prop1} and Proposition~\ref{prop2}$'$.
We have $\sigma_2(a) = \sigma(a) - \dsup[0]{a}$,
and ${\sigma_2}'(a_0)=\sigma'(a_0)-\dsup[0]{a_0}
=N+\sigma(a)-(N + \dsup[0]{a} ) = \sigma_2(a)$.
The number of chips in $a_1$ becomes:
${\sigma_2}'(a_1)=\sigma'(a_1)+N-\dsup[0]{a}=N+\sigma(a)-\dsup[0]{a}=\sigma_2(a)$.

With the same arguments, we obtain that any firing sequence $s$ of $C$ corresponds to a valid firing sequence of $C'$
satisfying Proposition~\ref{prop1} and Proposition~\ref{prop2} or \ref{prop2}$'$.
Notice that, for each configuration $\sigma$ of $C$ corresponding to a configuration $\sigma'$ in $C'$,
the vertices that can be fired in both configurations are the same.
Therefore, there is no other valid firing sequence in $C'$ than those corresponding to firing sequences of $C$.
Therefore \confspace{C'} is isomorphic to \confspace{C}.

\end{proof}

For simple MCFGs, the \shot{} of a configuration is defined in the same way as for CFGs.
Likewise, the configuration space of any convergent MCFG is isomorphic to the set of its \shot{}s ordered by inclusion.
We now show that any MCFG is equivalent to a CFG.

\begin{theoreme}
Any convergent MCFG  is equivalent
to a  CFG.
\end{theoreme}
\begin{proof}
Let $C$ be a convergent MCFG. From Theorem~\ref{cormutesimple} 
we can suppose that $C$ is simple.
Let $G=(V,E)$ and $\sigma$ be the initial
support graph and configuration of $C$. To
prove our result, we will use $C$ to construct an equivalent (simple) CFG.
Given a vertex $v$, we define the initial degree of $v$ to be its degree in $G$.
When $v$ is fired (if fired at all), a mutation occurs and its outdegree changes.
Since $C$ is simple, $v$ will be fired only once before the end of an
execution, so its outdegree changes only once.
We define $v$'s \emph{final} degree to be its degree in the final configuration
(in the case where $v$
is not fired during an execution, the initial and final degrees are the same). 

In the final configuration, $v$ contains less chips than its final outdegree (otherwise $v$ could be fired and we would not be in the final configuration).
Two cases can occur:
either $v$ contains also less chips than its \emph{initial} outdegree, or not.
In the first case, if we suppress the mutation of $v$
(the outgoing edges of $v$ are the same before and after the mutation), $v$ still cannot be fired in the final configuration.
If this is the case for all vertices, we can suppress all mutations to obtain a simple classical CFG.

In the other case, there exists a vertex $w$ which contains at the end of an execution more
chips than its initial outdegree (which
means that the outdegree of the vertex has been increased by the
mutation). 
We will modify $C$ to prevent this from happening.
We 
proceed in the following way: let \dplus{w} be the initial outdegree of $w$,
\conffin{w} be the number of chips in $w$ in the final configuration, and $n$
be the number of chips needed to fire $w$ (if $\sigma(v) \ge d$, then
$n=0$). 
We
add $\conffin{w}+1-\dplus{w}$ edges from $w$ to the sink, and
we set the initial
configuration of $w$ to $\conffin{w}+1-n$ chips. 
Figure~\ref{figcfgmut-cfg} shows an example of this construction.
After this modification, $w$ still needs $n$ chips to be fired, therefore the first times at which $w$ can be fired are not modified, and $w$ now 
contains less chips at the end of an execution than its initial outdegree.
We do this for every vertex that does not satisfy this property, obtaining thus a MCFG that can be considered as a CFG.
\end{proof}

\begin{figure}
\mbox{} \hfill
\begin{minipage}{20mm}
\begin{picture}(20,20)
\node(a)(5,15){3}
\node(b)(5,5){0}
\node[fillgray=0.85](puits)(15,10){}
\drawedge[curvedepth=2](a,b){} \drawedge[curvedepth=-2](a,b){}
\drawedge(a,b){}
\drawedge(b,puits){}
\end{picture}
\end{minipage}
\hfill \begin{minipage}{10mm}
 $\longrightarrow$ \end{minipage}
\begin{minipage}{20mm} 
\begin{picture}(20,20)
\node(a)(5,15){0}
\node(b)(5,5){2}
\node[fillgray=0.85](puits)(15,10){}
\drawedge(a,b){}
\drawedge(b,puits){}
\drawedge[curvedepth=2](b,puits){} \drawedge[curvedepth=-2](b,puits){}
\end{picture}
\end{minipage}
\hfill
\rule[-10mm]{.5mm}{20mm}
\hfill
\begin{minipage}{20mm} 
\begin{picture}(20,20)
\node(a)(5,15){3}
\node(b)(5,5){2}
\node[fillgray=0.85](puits)(15,10){}
\drawedge(a,b){}
\drawedge[curvedepth=2](a,b){} \drawedge[curvedepth=-2](a,b){}
\drawedge(b,puits){}
\drawedge[curvedepth=2](b,puits){} \drawedge[curvedepth=-2](b,puits){}
\end{picture}
\end{minipage}
\hfill \begin{minipage}{10mm} $\longrightarrow$ \end{minipage}
\begin{minipage}{20mm} \begin{picture}(20,20)
\node(a)(5,15){0}
\node(b)(5,5){2}
\node[fillgray=0.85](puits)(15,10){}
\drawedge(a,b){}
\drawedge[curvedepth=2](a,b){} \drawedge[curvedepth=-2](a,b){}
\drawedge(b,puits){}
\drawedge[curvedepth=2](b,puits){} \drawedge[curvedepth=-2](b,puits){}
\end{picture} \end{minipage}
\hfill \mbox{}
\caption{An example of a simple MCFG transformed into a CFG.
On the left, the initial and final configurations of the MCFG.
Notice that the lower vertex contains in the final configuration more chips than its initial outdegree.
On the right, the initial and final configuration of the corresponding CFG.}
\mylabel{figcfgmut-cfg}
\end{figure}

Although the MCFG is a generalization of the CFG, this theorem shows that in fact it generates no more lattices than the usual Chip Firing Game.
As we will see in the next section, the Chip Firing Game generates more lattices than the Abelian Sandpile Model, which is a CFG on an undirected graph.
Finding out if a model defined as a Mutating Abelian Sandpile Model generates the same lattices as the usual ASM would help understand these differences better.

\section{Comparison of CFGs and ASMs} \mylabel{seccfgasm}

In this section, we compare the classes of configuration spaces induced by CFGs and ASMs.
Since we know from Section~\ref{subsubdef} that any ASM is equivalent to a CFG,
we try to determine at which conditions a CFG is equivalent to an ASM.
We will show that this is always the case when the CFG is simple and its support graph has no cycle,
by giving an algorithm to transform such a CFG into an ASM.
Since we know that any distributive lattice is the configuration space of a simple CFG with no cycle~\cite{MPV01}, 
we obtain as a corollary that the class of lattices induced by ASM 
contains the distributive lattices.
However, we will prove that not all CFGs are equivalent to an ASM.
Since examples of lattices in \lasm{} that are not distributive appear often in the study of the model, we obtain the result that \lasm{}
is strictly between the distributive lattices and the lattices induced by CFG,
which is surprising because these two classes are very close to one another.

We will try to transform CFGs into ASMs using local transformations on the support graph of a CFG that do not change its configuration space.
Using a combination of these transformations, we will see that
it is possible in some cases to obtain a CFG such that its support graph contains one edge $(u,v)$ for each edge $(v,u)$, \ie{} the graph can be viewed as undirected,
 so that we have an ASM.
We begin by giving the two basic transformations we will use.
They do not change the vertex set of a CFG, only its edge set and its initial configuration.
Thanks to  Proposition~\ref{propcaracconfspace}, we will prove that these modifications preserve the configuration space by checking that the first times each vertex can be fired are not changed, and that the CFG remains simple
(indeed, in some cases, a modification does not change the first  times at which a given vertex can be fired, but allows this vertex to be fired more than once).

The first modification is a basic one, that will be used several times with small variations in the sequel.\\

\noindent
\textbf{Modification 1: Grounding}\\
The grounding modification applied with factor $n$ to a vertex $v$ of a CFG consists in
adding $n$ chips on $v$ in the initial configuration, and adding $n$ edges from $v$ to the sink.\\

\noindent
\textbf{Modification 2: Multipliying}\\
The multiplying modification consists in multiplying by an integer factor $n$ the indegree and initial configuration of a given vertex $v$, without modifying the rest of the CFG.
It consists in:
\begin{itemize}
\item multiplying by $n$ the initial configuration of $v$,
\item adding $(n-1)\cdot\dplus{v}$ edges from $v$ to the sink,
\item for each immediate predecessor $u$ of $v$, 
     adding $(n-1)\cdot d(u,v)$ edges $(u,v)$, and
     adding as much chips to the initial configuration of $u$.
\end{itemize}

The next two lemmas state that these modifications do not change the configuration space of the CFG to which they are applied.
The first result being obvious, its proof is ommitted.

\begin{lemme} \mylabel{lemgr}
The CFG obtained by applying the \gr{} to a CFG $C$ is equivalent to $C$.
\end{lemme}

\begin{lemme} \mylabel{lemmul}
The CFG obtained  by applying the \mul{} to a CFG $C$ is equivalent to $C$.
\end{lemme}
\begin{proof}
Let $C'$ be the CFG obtained from $C$ by applying the \mul{} to a vertex $v$ with a factor $n$.
\confspace{C'} is isomorphic to \confspace{C} since:
\begin{itemize}
\item The outdegree and the initial configuration of $v$ have been multiplied by $n$, thus $v$ needs $n$ times more chips to be fired in $C'$ than in $C$.
Since, for each immediate predecessor $u$ of $v$ the number of edges from $u$ to $v$ has also been multiplied by $n$, 
the first times at which $v$ can be fired do not change.
We can easily see that $v$ can be fired only once in $C'$, since:
$$\begin{array}{cl} 
\sigma'(v)+\sum d'(u,v) & = n\sigma(v) + n\sum d(u,v) \\
	&  < nd^+(v) \\
	& < {d'}^+(v)
\end{array}$$
\item the modifications made to the predecessors of $v$ (the modification of the initial configuration and the edges added towards $v$) 
are similar to applying the \gr{} to these vertices
(with $v$ playing the role of the sink).
\end{itemize}
\end{proof}

\begin{algorithm}
\In{A simple CFG $C$ with support graph $G=(V,E)$ and initial configuration $\sigma_0$, such that $G$ has no cycle}
\Out{An ASM equivalent to $C$}
Compute \confspace{C};\\
For each $v\in V\setminus\ens{\bot}$, 
$D[v]\gets\max_{\sigma\in\confspace{C}} (\dplus{v}-\sigma(v))$;\\
Mark $\bot$;\\
 $L\gets$~empty-list;\\
\While{there are unmarked vertices}{
  Choose an unmarked vertex $v$ with all successors marked;\\
  Mark $v$;\\
  $L \gets L,v$;
}
\textbf{Step 1:} \While{$L$ is not empty}{
  $v\gets$ \mbox{head}(L)\\
  \If{$D[v]\le d(v)$}{
    Apply the \mul{} to $v$ with  factor 
   \mbox{$\left\lceil(d(v)+1)/D[v]\right\rceil$};\\
    For each edge $(u,v)$ added by the \mul{}, add 1 to $D[u]$;}
  $L\gets \mbox{tail}(L)$;
  }
\textbf{Step 2:} \ForEach{ $v\in V\setminus\ens{\bot}$}{
  \If{$\confinit{v}+\dmoins{v}+d(v) \ge 2\dplus{v}$}{
    Apply the \gr{} to $v$ with  factor 
    \mbox{$\left\lceil(d(v) +\dmoins{v}+1)/
         (2\cdot\dplus{v}-\sigma_0(v))\right\rceil$};
  }
} 
\textbf{Step 3:} \ForEach{ $v\in V\setminus\ens{\bot}$}{
  \ForEach{edge $(u,v)$ in $E$}{
    Add one chip to the initial configuration of $v$;\\
    Add one edge $(v,u)$ to $E$;
  }
}
\caption{Transformation of a simple CFG with no cycle into an equivalent ASM}
\mylabel{algcfgasm}
\end{algorithm}

\noindent
Now we give the main theorem of this section:
\begin{theoreme} \mylabel{thmcnasm}
Algorithm~\ref{algcfgasm} transforms a simple CFG with no cycle into an equivalent ASM in linear time with respect to the size of its configuration space.
\end{theoreme}
\begin{proof}
The goal of the algorithm is to obtain an ASM from a CFG, by adding a \emph{reverse} edge to each edge in the graph
(\ie{} add an edge $(v,u)$ for each edge $(u,v)$).
We recall that, when a CFG is modified, we obtain an equivalent CFG if:
\begin{itemize}
\item the first times at which each vertex can be fired do not change, and
\item each vertex can be fired only once during an execution.
\end{itemize}
The purpose of the first two steps of the algorithm is to prepare the graph for this reversal of edges.
We first show that, at each step of the algorithm, the current CFG is equivalent to $C$.\\

\noindent
\textbf{Step 1}
In this step the indegree of a vertex is increased so that it exceeds by a certain factor its outdegree.
This is possible because the graph has no cycle, and the graph has only one sink.
In this step, the only modification applied is the \mul{},
therefore by Lemma~\ref{lemmul},
at each step the current CFG is equivalent to $C$.
At the end of this step, we have, for each vertex $v$ and for all configuration $\sigma$:
$$\dplus{v}-\sigma(v) > d(v).$$

\noindent
\textbf{Step 2}
The only modification applied in this step is the \gr{},
therefore, by Lemma~\ref{lemgr}, at each step the current CFG is equivalent to $C$.
At the end of this step we have, for each vertex $v$:
$$\confinit{v}+\dmoins{v}+d(v)<2\dplus{v}.$$

\noindent
\textbf{Step 3}
In this step an edge $(u,v)$ is added for each pre-existing edge $(v,u)$.
The resulting CFG is equivalent to $C$ because:
\begin{itemize}
\item a chip is added to the initial configuration of $u$ for each edge $(u,v)$ added. This is similar to applying the \gr{} to $u$, with $v$ playing the role of the sink.
Figure~\ref{figCN1}~(1-3) illustrates why this is necessary.
\item After Step~2 we know that $\confinit{v}+\dmoins{v}+d(v)<2\dplus{v}$.
      This means that, after we have added $d(v,u)$ edges from $u$ to $v$, we still have $\confinit{v}+\dmoins{v}<2\dplus{v}$, which means that $v$ can be fired only once in the resulting CFG.
Figure~\ref{figCN1} (4) illustrates this.
\item After Step~1 we have, for all configuration $\sigma$: $\dplus{v}-\sigma(v) > d(v)$.
      Therefore adding edges from $u$ to $v$ does not change the first times at which $v$ can be fired:
      there is no configuration $\sigma$ in which $v$ cannot be fired, and in which adding edges from $u$ to $v$ allows $v$ to be fired.
Figures~\ref{figCN3} and~\ref{figCN4} illustrate this.
\end{itemize}



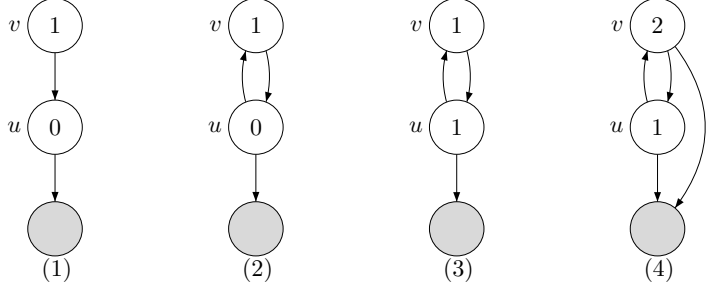
\begin{figure}
\begin{center}
\hfill
\scalebox{0.90}{
\begin{picture}(10,35)
\node(a)(0,35){1}
\node(b)(0,20){0}
\node[fillgray=0.85](c)(0,5){}
\gasset{ExtNL=y,NLangle=180,NLdist=1}
\nodelabel(a){$v$}
\nodelabel(b){$u$}
\drawedge(a,b){}
\drawedge(b,c){}
\put(-2,-2){(1)}
\end{picture}}
\hfill
\scalebox{0.90}{
\begin{picture}(10,35)
\node(a)(0,35){1}
\node(b)(0,20){0}
\gasset{ExtNL=y,NLangle=180,NLdist=1}
\nodelabel(a){$v$}
\nodelabel(b){$u$}
\node[fillgray=0.85](c)(0,5){}
\drawedge[curvedepth=2](a,b){}
\drawedge[curvedepth=2](b,a){}
\drawedge(b,c){}
\put(-2,-2){(2)}
\end{picture}}
\hfill
\scalebox{0.9}{
\begin{picture}(10,35)
\node(a)(0,35){1}
\node(b)(0,20){1}
\node[fillgray=0.85](c)(0,5){}
\gasset{ExtNL=y,NLangle=180,NLdist=1}
\nodelabel(a){$v$}
\nodelabel(b){$u$}
\drawedge[curvedepth=2](a,b){}
\drawedge[curvedepth=2](b,a){}
\drawedge(b,c){}
\put(-2,-2){(3)}
\end{picture}}
\hfill 
\scalebox{0.9}{
\begin{picture}(10,35)
\node(a)(0,35){2}
\node(b)(0,20){1}
\node[fillgray=0.85](c)(0,5){}
\gasset{ExtNL=y,NLangle=180,NLdist=1}
\nodelabel(a){$v$}
\nodelabel(b){$u$}
\drawedge[curvedepth=2](a,b){}
\drawedge[curvedepth=2](b,a){}
\drawedge(b,c){}
\drawedge[curvedepth=7](a,c){}
\put(-2,-2){(4)}
\end{picture}}
\hfill \mbox{}
\caption{\textbf{(2)}: Simply adding an edge $(u,v)$ to the CFG (1) changes the first times at which $u$ can be fired:  $u$ can not be fired at all.
\textbf{(3)}: therefore, we add a chip to the initial configuration of $u$, and $u$ can be fired after $v$, as in (1).
\textbf{(4)}: In~(3) $v$ can be fired more than once. Therefore,
 we apply the \gr{} to $v$ (\ie{} we add a chip to the initial configuration of $v$ and we add one edge $(v,\bot)$) before adding the edge $(u,v)$. The resulting CFG is equivalent to (1).}
\mylabel{figCN1}
\end{center}
\end{figure}



\begin{figure}
\hfill
\scalebox{0.9}{
\begin{picture}(15,35)
\node(a)(0,35){1}
\node(b)(15,35){1}
\node(v)(0,20){0}
\node(u)(15,10){1}
\node[fillgray=0.85](puits)(0,0){}
\gasset{ExtNL=y,NLdist=1,NLangle=180}
\nodelabel(a){$a$}
\nodelabel[NLangle=0](b){$b$}
\nodelabel(v){$v$}
\nodelabel[NLangle=0](u){$u$}
\drawedge(a,v){} \drawedge(b,u){} \drawedge(v,u){}
\drawedge[curvedepth=2](u,puits){} \drawedge[curvedepth=-2](u,puits){}
\end{picture}}
\hfill
\scalebox{0.85}{\input{figures/CN3.pstex_t}}
\hfill \mbox{}
\caption{A CFG and its configuration space. Notice that $v$ needs one chip to be fired, and that there exactly is one edge $(v,u)$.
After $b$ and $u$ have been fired
(the configuration is marked with a dot in the configuration space), only $a$ can be fired.
}
\mylabel{figCN3}
\end{figure} 

\begin{figure}
\hfill
\scalebox{0.9}{
\begin{picture}(15,35)
\node(a)(0,35){1}
\node(b)(15,35){1}
\node(v)(0,20){1}
\node(u)(15,10){2}
\node[fillgray=0.85](puits)(0,0){}
\gasset{ExtNL=y,NLdist=1,NLangle=180}
\nodelabel(a){$a$}
\nodelabel[NLangle=0](b){$b$}
\nodelabel(v){$v$}
\nodelabel[NLangle=0](u){$u$}
\drawedge(a,v){} \drawedge(b,u){} \drawedge[curvedepth=2](v,u){}
\drawedge[curvedepth=2](u,v){}
\drawedge[curvedepth=2](u,puits){} \drawedge[curvedepth=-2](u,puits){}
\drawedge(v,puits){}
\end{picture}}
\hfill
\scalebox{0.85}{\input{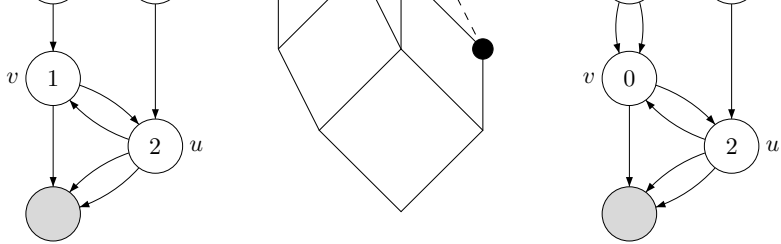}}
\hfill
\scalebox{0.9}{
\begin{picture}(15,35)
\node(a)(0,35){2}
\node(b)(15,35){1}
\node(v)(0,20){0}
\node(u)(15,10){2}
\node[fillgray=0.85](puits)(0,0){}
\gasset{ExtNL=y,NLdist=1,NLangle=180}
\nodelabel(a){$a$}
\nodelabel[NLangle=0](b){$b$}
\nodelabel(v){$v$}
\nodelabel[NLangle=0](u){$u$}
\drawedge[curvedepth=2](a,v){} \drawedge[curvedepth=-2](a,v){}
\drawedge(b,u){} \drawedge[curvedepth=2](v,u){}
\drawedge[curvedepth=2](u,v){}
\drawedge[curvedepth=2](u,puits){} \drawedge[curvedepth=-2](u,puits){}
\drawedge(v,puits){}
\end{picture}}
\hfill \mbox{}
\caption{If we want to add one edge $(u,v)$ to the CFG of Figure~\ref{figCN3},
we have first to
add one chip to the initial configuration of $u$, and to apply the \gr{} to $v$, as seen before.
We thus obtain the CFG shown on the left.
In this CFG,
it is possible to fire $v$ after $b$ and $u$,
and this adds the dotted lines to the configuration space.
Therefore, we must apply the \mul{} to $v$ to multiply its indegree by two, before we add the edge $(u,v)$. The resulting CFG (on the right) is equivalent to the original CFG (see Figure~\ref{figCN3}).}
\mylabel{figCN4}
\end{figure}

\end{proof}

\begin{corollaire}
Let $C$ be a simple CFG with support graph $G=(V,E)$ such that $G$ has no cycle.
Then $C$ is equivalent to an ASM.
\end{corollaire}

We have just shown that any simple CFG with no cycle is equivalent to an ASM.
This seems like a strong restriction, but in fact a large part of lattices of \lcfg{} are the configuration space of a simple CFG with no cycle.
Indeed, the next theorem states that this is the case for all distributive lattices, which are an important part of the ULD lattices, and therefore of \lcfg{}.

\begin{theoreme} \cite{MPV01}
Given any distributive lattice $L$, there exists a CFG $C$ with no cycle, such that \confspace{C} is isomorphic to $L$.
\end{theoreme}

\begin{corollaire}
\mylabel{corasmdist}
Any distributive lattice is the configuration space of an ASM.
\end{corollaire}

The distributive lattices are not the only ones that can be obtained 
both by CFG and ASM.
Indeed, most of the simple CFGs we have considered during our studies have no cycle, and are therefore equivalent to an ASM.
Moreover, we give in Figure~\ref{figCFGnonCNASM} an example of a CFG
which has a cycle in its support graph, but is equivalent to an ASM.

\begin{figure}
\mbox{} \hfill
\scalebox{0.8}{
\begin{picture}(20,40)
\node(a)(0,40){1}
\node(b)(20,40){1}
\node(c)(0,20){1}
\node(d)(20,20){1} \node[fillgray=0.85](puits)(10,0){}
\drawedge(a,c){} \drawedge(b,d){} \drawedge[curvedepth=3](c,d){}
 \drawedge[curvedepth=3](d,c){} \drawedge(c,puits){} \drawedge(d,puits){}
\end{picture}}
\hfill
\scalebox{0.8}{
\begin{picture}(20,40)
\gasset{AHnb=0}
\node(a)(0,40){2} \node(b)(20,40){2}
\node(c)(0,20){2} \node(d)(20,20){2} \node[fillgray=0.85](puits)(10,0){}
\drawedge(a,c){} \drawedge(b,d){} \drawedge(c,d){}
\drawedge(c,puits){} \drawedge(d,puits){}
\drawedge[curvedepth=-15](a,puits){}
\drawedge[curvedepth=15](b,puits){}
\end{picture}}
\hfill \mbox{}

\vskip -0.5cm
\mbox{} \hfill
\includegraphics[scale=0.7]{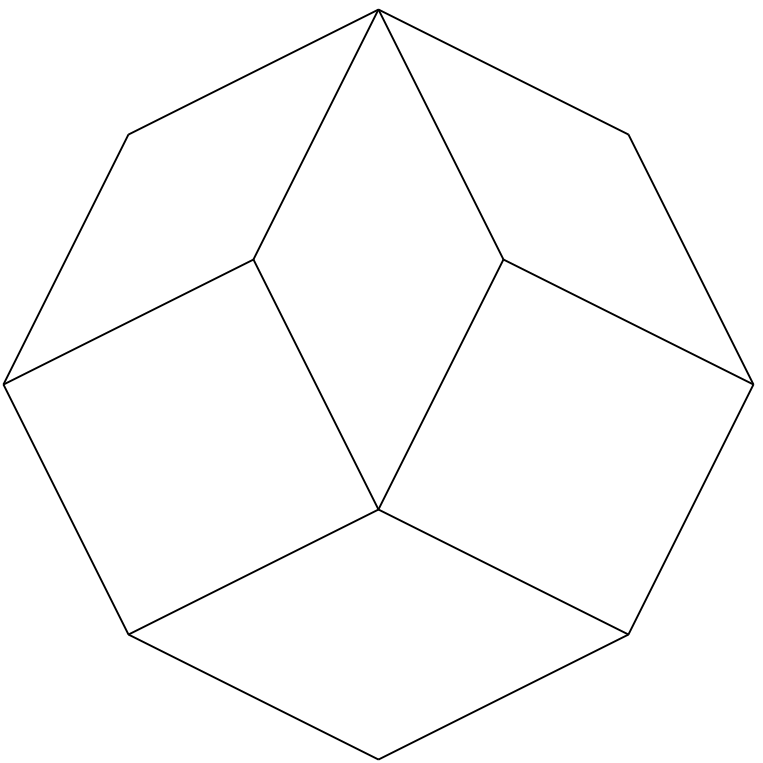}
\hfill \mbox{}
\caption{A CFG with a cycle in its support graph that is equivalent to an ASM.
We present here the CFG, an equivalent ASM, and their configuration space}
\mylabel{figCFGnonCNASM}
\end{figure}

However, we show that not all CFGs are equivalent to ASMs.

\begin{figure}
\begin{center}
\includegraphics[scale=0.5]{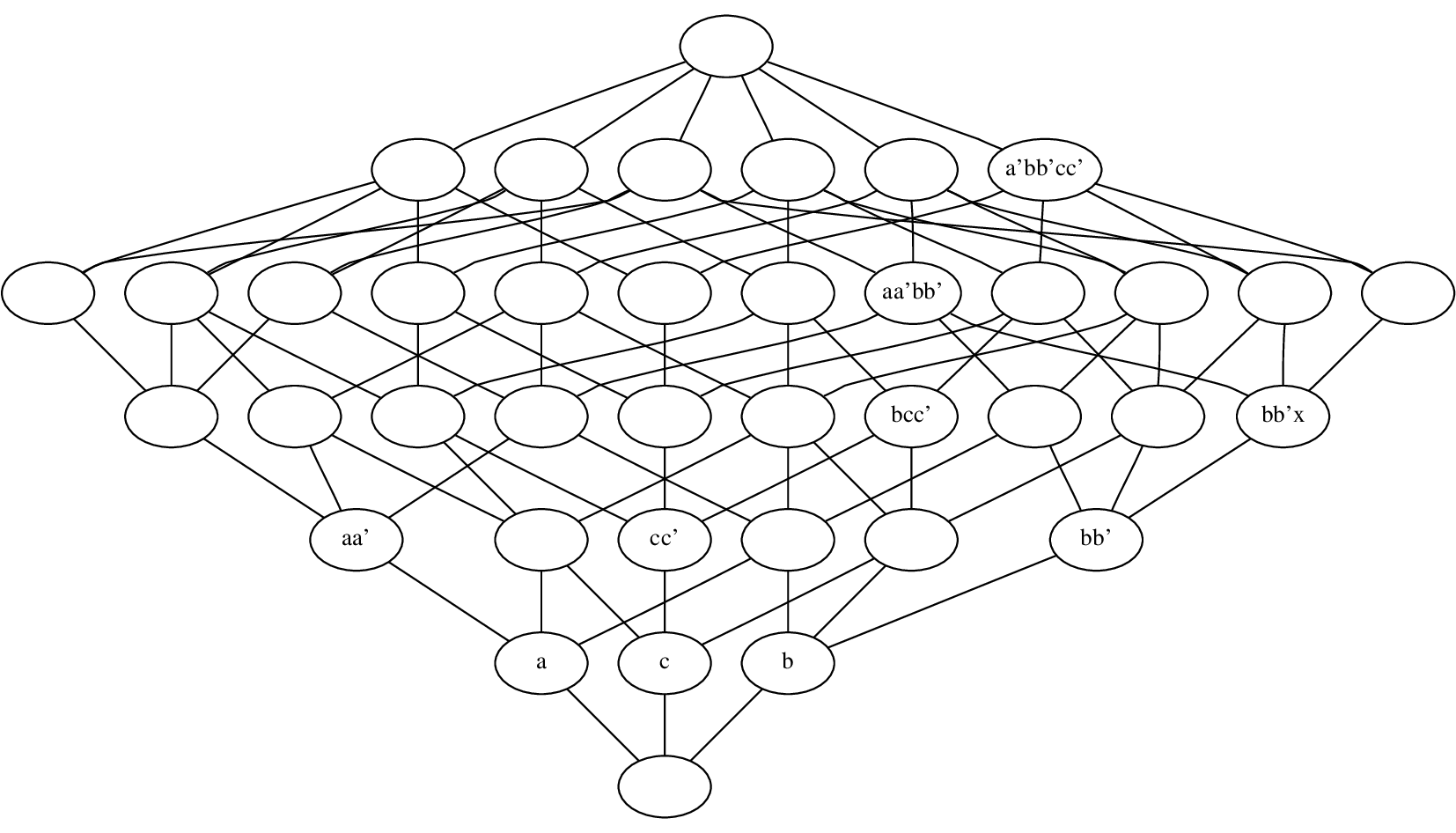}
\caption{A lattice that can be obtained by CFG but not by ASM}
\mylabel{figCFGnonASM}
\end{center}
\end{figure}

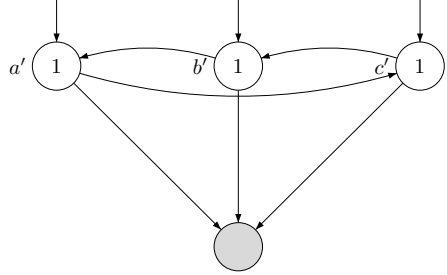
\begin{figure}
\begin{center}
\scalebox{0.8}{
\begin{picture}(60,60)
\node(a)(0,60){1}
\node(b)(30,60){1}
\node(c)(60,60){1}
\node(a')(0,30){1}
\node(b')(30,30){1}
\node(c')(60,30){1}
\gasset{ExtNL=y,NLangle=180,NLdist=1}
\nodelabel(b'){$b'$}
\nodelabel(a){$a$}
\nodelabel(a'){$a'$}
\nodelabel(c){$c$}
\nodelabel(b){$b$}
\nodelabel(c'){$c'$}
\node[fillgray=0.85](puits)(30,0){}
\drawedge(a,a'){}
\drawedge(b,b'){}  \drawedge(c,c'){}
\drawedge[curvedepth=-5](a',c'){}
\drawedge[curvedepth=-3](b',a'){} \drawedge[curvedepth=-3](c',b'){}
\drawedge(a',puits){} \drawedge(b',puits){}  \drawedge(c',puits){}
\end{picture}}
\caption{A CFG the configuration space of which is the lattice of Figure~\ref{figCFGnonASM}. Notice the cycle in the support graph.}
\mylabel{figCFGCFGnonASM}
\end{center}
\end{figure}

\begin{theoreme}
\mylabel{thmcfgpsasm}
$\lasm\varsubsetneq\lcfg.$
\end{theoreme}
\begin{proof}
We know that $\lasm\subseteq\lcfg$.
We will show that
the lattice $L$ of Figure \ref{figCFGnonASM}, which is the configuration space 
of the CFG of Figure \ref{figCFGCFGnonASM},  cannot be the configuration space
of any ASM
(\ie{} a CFG such that there exists an edge $(u,v)$ for each edge $(v,u)$).
This will show that the inclusion is strict.
The proof is in two steps: first we show that $L$ can be obtained only by a simple CFG 
(\ie{} a CFG such that its configuration space is $L$ \emph{must} have six different vertices, each fired once during an execution, and a sink);
we have to show that $L$ can only be obtained by a simple CFG, otherwise 
we must  prove also that $L$ cannot be the configuration space of an ASM with five vertices or less.
Then we prove some inequalities on the number of edges between two pairs of vertices, from which we will show that there exist two vertices 
$a$ and $b$ such that  the number of edges $(a,b)$ can never be the same as the number of edges $(b,a)$.
Therefore we can never consider that the support graph of the CFG is undirected, and this proves that $L$ cannot be obtained by ASM.

Let us prove that $L$ can be only obtained by a simple CFG.
The minimal element has three immediate successors, so there are three different vertices that can be fired in the initial configuration.
Let us call them $a,b$ and $c$.
The configuration denoted by $a$ has three immediate successors, so there are three different vertices that can be fired in configuration $a$.
Two of them are $b$ and $c$, and the third one, that we denote by $a'$, is therefore a vertex different from $b$ and $c$.
$a'$ is either a fourth vertex (that can be fired immediatly after $a$),
or $a$ itself (which would mean that $a$ can be fired twice in a row).
The same argument holds for $b$ and $c$.
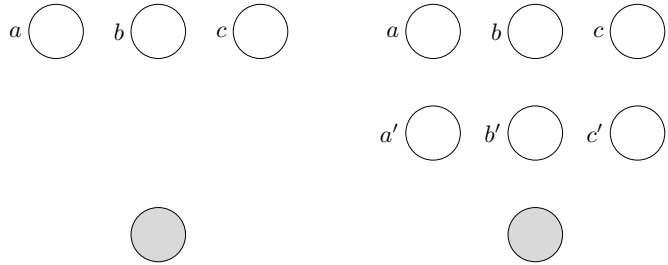
\begin{figure}
\begin{center}
\mbox{} \hfill
\scalebox{0.9}{
\begin{picture}(30,30)
\gasset{ExtNL=y,NLangle=180,NLdist=1}
\node(a)(0,30){$a$} \node(b)(15,30){$b$} \node(c)(30,30){$c$}
\node[fillgray=0.85](puits)(15,0){}
\end{picture}}
\hfill
\scalebox{0.9}{
\begin{picture}(30,30)
\gasset{ExtNL=y,NLangle=180,NLdist=1}
\node(a)(0,30){$a$} \node(b)(15,30){$b$} \node(c)(30,30){$c$}
\node[fillgray=0.85](puits)(15,0){}
\node(a')(0,15){$a'$} \node(b')(15,15){$b'$} \node(c')(30,15){$c'$}
\end{picture}}
\hfill \mbox{}
\caption{The possible vertex sets of a CFG that generates the lattice of Figure~\ref{figCFGnonASM}}
\mylabel{fig36sommets}
\end{center}
\end{figure}
So the vertex set of a CFG $C$ such that $L=\confspace{C}$ is something in-between
the two vertex sets of Figure \ref{fig36sommets}.

We will now prove that $a',b'$ and $c'$ are different from $a,b$ and $c$, and that they
are all distinct from each other.
First we prove that $a'$ is different from $a$:
the configuration $bb'$ has three immediate successors, therefore there are three different vertices that can be fired in configuration $bb'$.
Two of them are $a$ and $c$, so the third one (denoted by $x$, and the firing of which leads to the configuration labelled by $bb'x$) is either $a'$ (in which case $a'$ is distinct from $a$) or $c'$ (in which case $c'$ is distinct form $c$).
We observe that the join of $aa'$ and $bb'$ (labelled by $aa'bb'$) is greater than $bb'x$, which means that $\ens{bb'x}\subseteq\ens{aa'bb'}$.
From this we conclude that the third vertex is $a'$, and therefore $a'$ is distinct from $a$.
The same reasonning can be applied to $b'$ and $c'$.
Now we show that $a'$, $b'$ and $c'$ 
are distinct from each other:
if $a'$ and $b'$ are the same vertex, the join of $aa'$ and $bb'$ would be
reached from $aa'$ by the firing of $b$, and from $bb'$ by the firing of $a$.
Therefore, $aa'bb'$ would be
 at distance one of $aa'$ and of $bb'$, whereas it is at distance two.
Therefore, $a'$ and $b'$ are two different vertices.
With the same argument, we can show that $a'$, $b'$ and $c'$ are all distinct from each other, therefore
$L$ can only be obtained by a  CFG where six different vertices are fired during an execution.

Now we prove that $L$ cannot be the configuration space of any ASM.
First we claim that $d(b',a') > d(c',a')$.
Let $n$ be the number of chips that $a'$ needs to be fired.
We know that $a'$ can be fired after $b$ and $b'$. So we know that
$n\le d(b,a')+d(b',a')$.
Suppose our claim is not true, \ie{} $d(c',a') \ge d(b',a')$. 
Then  we have $n\le d(c',a')+d(b',a')$, and therefore
$n\le d(c',a')+d(b',a') + d(c,a')$. 
This means that $a'$ can be fired after $c,c'$ and $b'$.
But at point $b'cc'$, only two vertices can be fired, and they are $a$ and $b'$.
 Therefore we obtain a contradiction, and we must have $d(b',a') > d(c',a')$.
By similar arguments we can establish that $d(c',b') > d(a',b')$ and 
$d(a',c') > d(b',c')$.
Now if $L$ was the configuration space of an ASM, we would have
$d(a',b') = d(b',a')$,
$d(a',c') = d(c',a')$ and
$d(b',c') = d(c',b')$,
which is in contradiction with the inequalities established above.
We conclude that $L$ cannot be the configuration space of an ASM.
\end{proof}

By combining Theorem~\ref{thmcfgpsasm} and Corollary~\ref{corasmdist},
we obtain that \lasm{} is situated strictly between the class of distributive lattices and \lcfg{}.
This shows the complexity of the problems raised by the Chip Firing Game and the Abelian Sandpile Model in lattice theory:
\lcfg{} and \lasm{} are both between the distributive and ULD lattices,
while there is no previously known lattices class satisfying this condition.

\section*{Conclusion and perspectives}
In this paper, we have studied from the lattice point of view
the  Chip Firing Game and two
 closely related models: the Abelian Sandpile Model and the Mutating Chip Firing Game.
It was already known that all these models generate lattices, and that these lattices are in the class of Upper Locally Distributive (ULD) lattices.
Our goal was to 
characterize the classes of lattices \lcfg{}, \lasm{} and \lmcfg{} induced by each model (\ie{}
determine, given a ULD lattice, by which model(s) it can be obtained).
We have first shown that the mutating CFG and the CFG generate exactly the same set of lattices,
by giving a way to transform a MCFG into a classical CFG.
This means that, although the MCFG is a generalization of the CFG,
it brings no more information in lattice terms.

We know that every ASM is equivalent to a CFG,
and we have given an example of a CFG that is not equivalent to any ASM.
This implies that $\lasm{}\varsubsetneq \lcfg{}$.
Then we have given a sufficient but not necessary condition at which a CFG can be transformed into an ASM, 
from which it can be shown that the class $D$ of distributive lattices is included in \lasm{}.
The class of lattices induced by CFG being somewhere between the distributive and the ULD lattices,
we obtain that the class of lattices induced by ASM is a new class between
the distributive lattices and
 \lcfg{}.
In other words, we have proved the following relation:
$$D\varsubsetneq \mbox{\lasm{}}\varsubsetneq \mbox{\lcfg{}}=\mbox{\lmcfg{}} \varsubsetneq \mbox{ULD}$$

This illustrates the complexity of the problem in lattice theory, because there exists no known lattice class between the distributive and ULD lattices.
Therefore, an exact characterization of the two classes \lasm{} and \lcfg{} would be a very interesting result in lattice theory.

The CFG and the ASM share the same definition, except that the first is defined on a directed graph, and the latter on an undirected graph.
This might seem like a strong difference, and the first idea that comes to mind is that \lasm{} is much smaller than \lcfg{}.
The sufficient but not necessary condition we have given at which a CFG is equivalent to an ASM is that the support graph must contain no directed cycle.
This is a strong restriction, but it nonetheless implies that \lasm{} is a very significant part of \lcfg{} 
(it contains more than the distributive lattices class, which is a large part of the ULD class).

The fact that three important models, used in various domains like physics, computer science and social science, all induce strongly structured sets precisely situated between two classical types of lattices, shows the importance of the use of order theory in the context of dynamical models studies.
In our case it is even more interesting to notice that the models introduce new classes of lattices which one may study from the order theoretical point of view.

We have seen that the difference between these models does not lie in the fact that the graph is directed or not, but on the existence or not of directed cycles in the graph.
Simple CFGs with no cycle in their support graph are very strongly linked with the class $\cal{CN}$ of lattices obtained by a sequence of duplications of convex sets.
It would be a natural sequel of this work to find out in all lattices of $\cal{CN}$ induced by CFG can be obtained by ASM.
However, if this is true this does not give the result that $\lasm{}=\lcfg{}\cap\cal{CN}$,
because we have examples of lattices that are in $\lasm{}$, but not in $\cal{CN}$.

Finally, lattices also appear in some other kinds of discrete dynamical models defined in the context of tiling theory:
for some classes of tiling problems, one can define a local rearrangement of tiles, called \emph{flip}, 
which transforms a tiling of a given region into another tiling of the same region.
In some cases (mainly tilings with dominoes or with three lozenges \cite{Rem99b,BL01}),
it has been proved that the flip relation gives the distributive lattice structure to the set of all possible tilings of a given region.
A work similar to the one done in this paper, comparing the classes of lattices induced by each tiling problem, may lead to interesting result, both in tiling theory and in the study of discrete dynamical models.

\bibliographystyle{alpha}
\bibliography{biblio}
\end{document}